\input amstex
\documentstyle{amsppt}
\NoBlackBoxes

\document
\topmatter
    \title
        Conditioning by rare sources
    \endtitle
    \author
        M. Grendar
    \endauthor

    \address
        Department of Mathematics;
        FPV UMB;
        Tajovskeho 40;
        SK-974 01 Banska Bystrica;
        Slovakia.
        Institute of Mathematics and Computer Science;
        Banska Bystrica;
        Slovakia.
        Institute of Measurement Science;
        Bratislava;
        Slovakia
    \endaddress
    \email
        marian.grendar\@savba.sk
    \endemail
    \abstract In this paper we study the exponential decay of  posterior probability of a set of sources
    and conditioning by rare sources for both  uniform and general
    prior distributions of sources. The decay rate is determined by $L$-divergence
    and rare sources from a convex, closed set  asymptotically conditionally
    concentrate on an $L$-projection. $L$-projection on a linear family of sources belongs
    to $\varLambda$-family of distributions.  The results parallel those of
    Large Deviations for Empirical Measures (Sanov's Theorem and Conditional Limit Theorem).
    \endabstract
    \dedicatory
   To Mar, in memoriam. \\ 
   To George Judge,  on the occasion of his eightieth birthday.
    \enddedicatory
    \thanks
    I have  greatly benefited from valuable discussions with George Judge.
    Without implicating him, I am indebted to \v Lubom\'{\i}r Snoha for clarifications of a
    couple of important technical questions. Substantive comments and suggestions from two anonymous referees are gratefully acknowledged. Supported by  VEGA 1/7295/20 grant.  Date: August 31,
    2005. Revised: October 23, 2005. Second revision: November 25, 2005.
    \endthanks

\endtopmatter

\head 1. Introduction
\endhead

Information divergence minimization, which is also known as Relative
Entropy Maximization or MaxEnt method, has --  thanks to Large
Deviations Theorems for Empirical Measures -- gained a firm
probabilistic footing, which justifies its application in the area
of the convex Boltzmann Jaynes Inverse Problem (the
$\alpha$-problem, for short). For the $\beta$-problem -- an
'antipode' of the $\alpha$-problem -- Large Deviations Theorems for
Sources, which are presented here, single out the $L$-divergence
maximization method.

The paper is organized as follows: First, necessary terminology and
notation are introduced. A brief survey of Large Deviations Theorems
for Empirical Measures that includes Sanov's Theorem and a
Conditional Limit Theorem is given next. Then, a set-up for a study
of conditioning by rare sources is formulated and Sanov's Theorem
and the Conditional Limit Theorem for Sources are stated; under
various assumptions. Next, Theorems are proven for the continuous
case and the results are  applied to a criterion choice problem
associated with the $\beta$-problem. An End-Notes section points  to
relevant literature, mentions a motivation for the present work and
contains further discussion.

\head 2. Terminology and notation
\endhead

Let $\Cal{P}(\Cal{X})$ be a set of all probability mass functions
on a finite alphabet $\Cal{X} \triangleq \{x_1, x_2, \dots, x_m\}$
of $m$ letters. The support of $p \in \Cal{P}(\Cal{X})$ is a set
$S(p) \triangleq \{x: p(x) > 0\}$.

A probability mass function (pmf) from $\Cal{P}(\Cal{X})$ is
rational if it belongs to the set $\Cal{R} \triangleq
\Cal{P}(\Cal{X}) \cap \Bbb{Q}^m$. A rational pmf is $n$-rational, if
denominators of all its $m$ elements are $n$. The set of all
$n$-rational pmf's will be denoted by $\Cal{R}_n$.

Let $x_1, x_2, \dots, x_n$ be a sequence of $n$ letters, that is
identically and independently drawn from a source $q \in
\Cal{P}(\Cal{X})$. Type and $n$-type  are other names  for empirical
measures induced by a sequence of the length $n$. Formally, type
$\nu^n \triangleq [n_1, n_2, \dots, n_m]/n$, where $n_i$ is the
number of occurrences of $i$-th letter of the alphabet in the
sequence. Note that there are $\Gamma(\nu^n) \triangleq n!(
\prod_{i=1}^m n_i!)^{-1}$ different sequences of length $n$, which
induce the same type $\nu^n$. $\Gamma(\nu^n)$ is called the
multiplicity of type. Finally, observe that $\nu^n$ is $n$-rational;
$\nu^n \in \Cal{R}_n$.

Let $\Cal{\Pi}, \Cal{Q} \subseteq \Cal{P}(\Cal{X})$.
$\Cal{\Pi}_n \triangleq \Cal{\Pi} \cap \Cal{R}_n$ and $\Cal{Q}_n
\triangleq \Cal{Q} \cap \Cal{R}_n$. The former will be called set of
$n$-types $\nu^n$, the latter set of $n$-sources $q^n$.

The information divergence ($\pm$-relative entropy, Kullback-Leibler
distance etc.) $I(p||q)$ of $p$ with respect to  $q$ (both from
$\Cal{P}(\Cal{X})$) is 
$I(p||q) \triangleq \sum_\Cal{X} p \log\frac{p}{q}$, with
conventions that $0 \log 0 = 0$, $\log b/0 = + \infty$. The
information projection $\hat{p}$ of $q$ on $\Cal{\Pi}$ is  $\hat p
\triangleq \arg \inf_{p \in \Cal\Pi} I(p||q)$. The value of the
$I$-divergence at an $I$-projection of $q$ on $\Cal\Pi$ is denoted
by $I(\Pi||q)$.

On $\Cal P(\Cal X)$ topology induced by the standard topology on
$\Bbb{R}^m$ is assumed.

The support $S(\Cal C)$ of a convex set $\Cal C \subset
\Cal{P}(\Cal{X})$ is just the support of the member of $\Cal C$ for
which $S(\cdot)$ contains the support of  any other member of the
set.

The following families of distributions will be needed:

1) Linear family  $\Cal L(u, a) \triangleq \{p: \sum_\Cal{X} p(x)
u_j(x) = a_j, j = 1, 2, \dots, k\}$, where $u_j$ is a real-valued
function on $\Cal X$ and $a_j \in \Bbb{R}$.

2) Exponential family  $\Cal{E}(\rho, u, \theta) \triangleq \{p:
p(x) = z \rho(x) \exp({\sum_{j=1}^k \theta_j u_j(x)}), x \in
\Cal{X}\}$, where a normalizing factor $z \triangleq \sum_\Cal{X}
\rho(x) \exp({\sum_{j=1}^k \theta_j u_j(x)})$  and $\rho$ belongs to
$\Cal{P}(\Cal{X})$; $\theta_j \in \Bbb{R}$.

3) $\varLambda$-family  $\Cal{\varLambda}(\rho, u, \theta, a)
\triangleq \{p: p(x) = \rho(x)[1 - \sum_{j=1}^k \theta_k (u_j(x) -
a_j)]^{-1}, x \in \Cal{X}\}$.

The definitions of the families can be extended to continuous $\Cal
X$ in a straightforward way.

\smallskip

In what follows, $r \in \Cal{P}(\Cal{X})$ will be the 'true' source
of sequences and hence types.

\head 3. Conditioning by rare types
\endhead

It is convenient  to begin with a  brief survey of the Large
Deviations Theorems for Empirical Measures (Sanov's Theorem and a
Conditional Limit Theorem).

First, it is necessary to introduce the probability $\pi(\nu^n;
r)$ that the source $r$ generates an $n$-type $\nu^n$. The
probability that $r$ generates a sequence of $n$ letters $x_1,
x_2, \dots, x_n$ which induces a type $\nu^n$ is $\prod_{i=1}^m
(r_i)^{n\nu_i^n}$. As it was already mentioned, there is a number
$\Gamma(\nu^n)$ of sequences of length $n$, which induce the same
type $\nu^n$. The probability $\pi(\nu^n; r)$ that $r$ generates
type $\nu^n$ is thus $\pi(\nu^n; r) \triangleq \Gamma(\nu^n)
\prod_{i=1}^m (r_i)^{n\nu_i^n}$. Consequently, for $A \subseteq B
\subseteq \Cal P(\Cal X)$, $\pi(\nu^n \in A| \nu^n \in B; r) =
\frac{\pi(\nu^n \in A; r)}{\pi(\nu^n \in B;r)}$; provided that
$\pi(\nu^n \in B;r) \neq 0$.

$\Cal\Pi$ is rare if it does not contain $r$. Given that the source
$r$ produced an $n$-type from rare $\Cal\Pi$, it is of interest to
know how the conditional probability/measure spreads among the rare
$n$-types from $\Pi$; especially as $n$ grows beyond any limit. For
the rare set of a particular form, this issue is answered by
Conditional Limit Theorem (CoLT) which is also known as Conditional
Weak Law of Large Numbers.

CoLT can be established by means of Sanov's Theorem (ST).

\proclaim{ST} (\cite{6} Thm 3)   Let $\Cal\Pi$ be a set such that
its closure is equal to the closure of its interior. Let $r$ be
such that $S(r) = \Cal X$. Then,
$$
\lim_{n\rightarrow\infty} \frac{1}{n}\log\pi(\nu^n \in \Pi; r) = -
I(\Pi|| r).
$$
\endproclaim

Sanov's Theorem (ST) states that the probability $\pi(\nu^n \in
\Cal\Pi; r)$ decays exponentially fast, with the decay rate given by
the value of the information divergence at an $I$-projection of the
source $r$ on $\Cal\Pi$.

\proclaim {CoLT} (\cite{8} Thm 4.1, \cite{2} Thm 12.6.2)
Let $\Cal\Pi$ be a convex, closed rare set. 
Let $B(\hat p, \epsilon)$ be a closed $\epsilon$-ball defined by the
total variation metric, centered at $I$-projection $\hat p$ of $r$
on $\Cal\Pi$. Then for any $\epsilon > 0$,
$$
\lim_{n \rightarrow \infty} \pi(\nu^n \in B(\hat p, \epsilon) \,|\,
\nu^n \in \Cal\Pi; r) = 1.
$$
\endproclaim

Informally, CoLT states that if a dense rare set admits a unique
$I$-projection, then asymptotically types conditionally concentrate
just on it. Thus, provided that for sufficiently large $n$ a type
from rare $\Cal\Pi$ occurred, with probability close to $1$ it is
just a type close to $\hat p$. Numeric examples of ST and CoLT can
be found at \cite{2}.

This suggests that, conditionally upon the rare $\Cal\Pi$, it is the
$I$-projection $\hat p$ rather than $r$, which should be considered
as the true {\it iid} source of data. Gibbs' Conditioning Principle
(GCP) - an important strengthening of CoLT - captures this
'intuition'; cf \cite{3}, \cite{7}.

If $S(\Cal L) = \Cal X$ then the $I$-projection $\hat p$ of $r$ on
$\Cal\Pi \equiv \Cal L$ is unique and belongs to the exponential
family of distributions $\Cal{E}(r, u, \theta)$; i.e., $\Cal{L}(u,
a) \cap \Cal{E}(r, u, \theta) = \{\hat p\}$.

\head 4. Conditioning by rare sources
\endhead

In the above setting there is  a fixed source $r$ and a rare set
$\Pi_n$ of $n$-types. We now consider a setting where the $n$-type
is unique, and there is a set $\Cal{Q}_n$ of rare $n$-sources of
the type.

Furthermore, $n$-sources $q^n$ are assumed to have prior
distribution $\pi(q^n)$. If from $\Cal{R}_n$  $n$-source $q^n$
occurs, then the source generates $n$-type $\nu^n$ with the
probability $\pi(\nu^n|q^n) \triangleq \Gamma(\nu^n) \prod_{i=1}^m
(q^n_i)^{n \nu^n_i}$.

We are interested in the asymptotic behavior of  the probability
$\pi(q^n \in B\,|\,(q^n \in \Cal Q)\wedge \nu^n)$ that if the
$n$-type $\nu^n$ and an $n$-source $q^n$ from a rare set $\Cal Q$
occurred, then the $n$-source belongs to a subset $B$ of $\Cal Q$.
Note that $\pi(q^n \in B\,|\,(q^n \in \Cal Q)\wedge \nu^n) =
\frac{\pi(q^n \in B|\nu^n)}{\pi(q^n \in \Cal Q|\nu^n)}$; provided
that $\pi(q^n \in \Cal Q|\nu^n) > 0$. The posterior probability
$\pi(q^n| \nu^n)$ is related to  the defined probabilities
$\pi(\nu^n|q^n)$ and $\pi(q^n)$ via Bayes's Theorem.

Asymptotic investigations will be first carried on  under the
assumption of uniform prior distribution of $n$-sources (Sect.
4.1). The assumption will be relaxed in Section 4.2. Within each
of the sections, two cases of convergence will be considered: a
static and a dynamic case. For the static case asymptotic
investigations are  carried over a subsequence of types, which are
$k$-equivalent to $\nu^{n_0}$. A type $\nu^{kn_0} \triangleq
[kn_1, \dots, kn_m]/kn_0$, $k \in \Bbb{N}$, is called
$k$-equivalent to $\nu^{n_0}$. The dynamic case assumes that there
is a sequence of $n$-types which converges in the total variation
to some $p \in \Cal{P}(\Cal{X})$. For each case what is meant by
rare source will be defined separately.

For $p, q \in \Cal{P}(\Cal{X})$, the $L$-divergence $L(q||p)$ of
$q$ with respect to $p$ is the map $L: \Cal{P}(\Cal X) \times
\Cal{P}(\Cal X) \rightarrow \Bbb R\cup \{-\infty\}$, $L(q||p)
\triangleq \sum_\Cal{X} p \log q$. The $L$-projection $\hat q$ of
$p$ on  set of sources $\Cal Q$ is: $\hat q \triangleq \arg
\sup_{q \in \Cal{Q}} L(q||p)$. The value of $L$-divergence at an
$L$-projection (i.e., $\sup_{q \in \Cal Q} L(q||p)$) is denoted by
$L(\Cal Q||p)$.

 \subhead 4.1 Uniform prior
 \endsubhead

Within this section it is assumed that $n$-sources have a uniform
prior distribution. Since there is total  $N = {n + m - 1 \choose
m-1}$ $n$-sources (cf. \cite{4}), the uniform prior probability
$\pi(q_n) = 1/N$, for all $q^n \in \Cal{R}_n$.

 \subsubhead 4.1.1 Static case
 \endsubsubhead

Let there be an $n_0$-type $\nu^{n_0}$. A set $\Cal Q$ of sources
is  rare if it does not contain $\nu^{n_0}$.

Sanov's Theorem for Sources (abbreviated $L$ST)
 is a
counterpart of the Sanov's Theorem for Types.
\proclaim {Static $L$ST}
 Let $\nu^{n_0}$ be a type. Let $\Cal Q$ be an open
set of sources. Then, for $n \rightarrow \infty$ over a subsequence
$n = kn_0$, $k \in \Bbb N$,
 $$
 \frac{1}{n} \log \pi(q^n \in \Cal Q| \nu^n) = L(\Cal Q||\nu^{n_0}) -
 L(\Cal P||\nu^{n_0}).
 $$
\endproclaim

\demo {Proof} Under the assumption of uniform  prior distribution of
of $n$-sources
$$
\log \pi(q^n \in \Cal Q| \nu^n) = 
\log\sum_{q^n \in \Cal Q} \prod_\Cal{X} (q^n)^{n \nu^n} -
\log{\sum_{q^n \in \Cal{P}} \prod_\Cal{X} (q^n)^{n\nu^n}}.
$$
Since $N < (n+1)^m$ (cf. Lemma 2.1.2 of \cite{7}), $\frac{1}{n_0}
\log \pi(q^{n_0} \in \Cal Q| \nu^{n_0})$ can be bounded from above
and below as:
$$\multline
L(\Cal{Q}_{n_0}||\nu^{n_0}) - L(\Cal{R}_{n_0}||\nu^{n_0}) -
\frac{m}{n_0} \log(n_0 +
1) \le \frac{1}{n_0} \log \pi(q^{n_0} \in \Cal Q| \nu^{n_0}) \le \\
\le L(\Cal{Q}_{n_0}||\nu^{n_0}) - L(\Cal{R}_{n_0}||\nu^{n_0}) +
\frac{m}{n_0} \log(n_0 + 1).
\endmultline
$$

Fix $p \in \Cal{P}(\Cal X)$. Equip $\Bbb R \cup \{-\infty\}$ with
the standard topology (i.e., the topology induced by the total
order). As for each open subset $A$ of $\Bbb R \cup \{-\infty\}$,
$L^{-1}(A)$ is an open subset of $\Cal{P}(\Cal X)$, the
$L$-divergence is continuous in $q$.

$\Cal Q$ is open by the assumption.

Thus, $L(\Cal{Q}_{n_0}||\nu^{n_0})$ converges to $L(\Cal
Q||\nu^{n_0})$ as $n \rightarrow \infty$, $n = k n_0$, $k \in \Bbb
N$. Also, $L(\Cal{R}_{n_0}||\nu^{n_0})$ converges to $L(\Cal
P||\nu^{n_0})$ for $n \rightarrow \infty$, $n = k n_0$, $k \in \Bbb
N$. \qed
\enddemo

The Law of Large Numbers for Sources ($L$LLN) is a direct
consequence of $L$ST.
 \proclaim {Static $L$LLN}
 Let $\nu^{n_0}$ be a type. Let $\hat q$ be $L$-projection of
 $\nu^{n_0}$ on $\Cal P(\Cal X)$. And let $B(\hat q, \epsilon)$ be a closed
 $\epsilon$-ball defined by the total variation metric, centered at $\hat
 q$. Then, for $\epsilon > 0$ and $n \rightarrow \infty$ over 
 the types which are $k$-equivalent with $\nu^{n_0}$,
 $$
 \pi(q^n \in {B}(\hat q, \epsilon) | \nu^n) = 1.
 $$
\endproclaim

\demo {Proof} Let $B^C(\hat q, \epsilon) \triangleq \Cal P(\Cal
X)\backslash B(\hat q, \epsilon)$. Since $B^C(\hat q, \epsilon)$ is
open by the assumption, $L$ST can be applied to it. Since $B^C
\subset \Cal P$, $L(B^C||\nu^{n_0}) - L(\Cal P||\nu^{n_0}) < 0$.
Thus, $\pi(q^n \in B^C(\hat q, \epsilon) | \nu^n)$ converges to $0$,
as $n \rightarrow \infty$ over a subsequence of $n = kn_0$, $k \in
\Bbb N$. $\square$
\enddemo

Obviously, the $L$-projection $\hat q$ of $\nu^{n_0}$ on $\Cal
P(\Cal X)$ is $\hat q \equiv \nu^{n_0}$.

$L$LLN is a special, unconditional case of the Conditional Limit
Theorem for Sources ($L$CoLT), which is  a consequence of $L$ST, as
well.

\proclaim {Static $L$CoLT} Let $\nu^{n_0}$ be a type. Let $\Cal Q$
be a convex, closed rare set of sources. Let $\hat q$ be the
$L$-projection of $\nu^{n_0}$ on $\Cal Q$ and let $B(\hat q,
\epsilon)$ be a closed $\epsilon$-ball defined by the total
variation metric, centered at $\hat q$. Then, for $\epsilon > 0$ and
$n \rightarrow \infty$ over a subsequence $n = kn_0$, $k \in \Bbb
N$,
 $$
 \pi(q^n \in {B}(\hat q, \epsilon) \, | \, (q^n \in \Cal Q) \wedge \nu^n) = 1.
 $$
\endproclaim

\demo {Proof} Let $B^C(\hat q, \epsilon) \triangleq \Cal P(\Cal
X)\backslash B(\hat q, \epsilon)$. Clearly,
$$
\log\pi(q^{n_0} \in B^C(\hat q, \epsilon) \, | \, (q^{n_0} \in \Cal
Q) \wedge \nu^{n_0}) = 
\log{\pi(q^{n_0} \in B^C|\nu^{n_0})} - \log{\pi(q^{n_0} \in
 \Cal Q|\nu^{n_0})}.
$$
Since both $B^C(\hat q, \epsilon)$ and $\Cal Q$ are
 open, $L$ST can be applied. As $B^C(\hat q, \epsilon) \subset \Cal
 Q$, $L(B^C||\nu^{n_0}) - L(\Cal Q||\nu^{n_0}) < 0$. Hence $\pi(q^n \in B^C|(q^n \in \Cal Q) \wedge \nu^n)$
converges to $0$, as $n \rightarrow \infty$ over a subsequence of $n
= kn_0$, $k \in \Bbb N$. Since under the assumptions on $\Cal Q$ the
$L$-projection of $\nu^{n_0}$ on $\Cal Q$ is unique, the claim of
the Theorem follows.
 $\square$
\enddemo

\example {Example} Let $\Cal X = \{1, 2, 3, 4\}$. Let $\Cal Q = \{q:
\sum_{x \in \Cal X} q(x) x = 1.7\}$. Let $n_0 = 10$ and $\nu^{n_0} =
[1, 1, 1, 7]/10$. The $L$-projection of $\nu^{n_0}$ on $\Cal Q$ is
$\hat q = [0.705,\ 0.073,\ 0.039,\ 0.183]$. Let $\epsilon = 0.1$.
The concentration of $n$-sources on the $L$-projection, which is
captured by the Static $L$CoLT, is for types $k$-equivalent to
$\nu^{n_0}$ ($k = 5, 10, 20, 30$) illustrated in Table 1.

 \midinsert
 \captionwidth{10cm}
 \topcaption{Table 1} Values of $\pi(q^n \in B(\hat q, \epsilon)|(q^n \in \Cal Q)\wedge
  \nu^n)$ for $n = k n_0$, $k = 5, 10, 20, 30$.
 \endcaption
$$
\vbox{
 \offinterlineskip
 \halign{
 \strut\vrule\qquad #\qquad & \vrule\qquad #\qquad \vrule \cr
 \noalign{\hrule}
  $n$ & $\pi(\cdot|\cdot)$ \cr
 \noalign{\hrule}
 \noalign{\hrule}
  50 & 0.868 \cr
 \noalign{\hrule}
 100 & 0.948 \cr
 \noalign{\hrule}
 200 & 0.994 \cr
 \noalign{\hrule}
 300 & 0.999 \cr
 \noalign{\hrule}
 }
}
$$
\endinsert
\endexample

The $L$-projection at the above Example can be found by means of the
following Proposition.

\proclaim {Proposition} Let $\Cal Q \equiv \Cal{L}(u, a)$. Let $p
\in \Cal{P}(\Cal X)$ be such that $S(p) = S(\Cal L)$.  Then the
$L$-projection $\hat q$ of $p$ on $\Cal Q$ is unique and belongs to
$\Cal\varLambda(p, u, \theta, a)$ family; i.e., $\Cal{L}(u, a) \cap
\Cal{\varLambda}(p, u, \theta, a) = \{\hat q\}$.
\endproclaim

\demo {Proof} In light of Theorem 9 of \cite{6}  it  suffices to
check that $\hat q = p [1 - \sum_{j=1}^k \theta_k (u_j(x) -
a_j)]^{-1}$, with $\theta$ such that $\hat q \in \Cal{L}(u, a)$,
satisfies:
$$
\sum_{S(p)} p \left(1 - \frac{q'}{\hat q}\right) = 0,
$$
for all $q' \in \Cal Q$, which is indeed the case. $\square$
\enddemo

 \subsubhead 4.1.2 Dynamic case
 \endsubsubhead

Let there be a sequence of $n$-types which converges in the total
variation to a pmf $p \in \Cal{P}(\Cal{X})$, denoted as $\nu^n
\rightarrow p$. In this case, a set $\Cal Q$ of sources is
 rare if it does not contain $p$.

\proclaim {Dynamic $L$ST}
 Let $\nu^n \rightarrow p$. Let $\Cal Q$ be an open
set of sources. Then,
 $$
\lim_{n \rightarrow \infty} \frac{1}{n} \log \pi(q^n \in Q| \nu^n) =
L(\Cal Q||p) - L(\Cal P||p).
 $$
\endproclaim

 \proclaim {Dynamic $L$LLN}
 Let $\nu^n \rightarrow p$. Let $\hat q$ be $L$-projection of
 $p$ on $\Cal P$. And let $B(\hat q, \epsilon)$ be a closed
 $\epsilon$-ball defined by the total variation metric, centered at $\hat
 q$. Then, for $\epsilon > 0$,
 $$
\lim_{n \rightarrow \infty} \pi(q^n \in {B}(\hat q, \epsilon) |
\nu^n) = 1.
 $$
\endproclaim

\proclaim {Dynamic $L$CoLT} Let $\nu^n \rightarrow p$. Let $\Cal Q$
be a convex, closed rare set of sources. Let $\hat q$ be the
$L$-projection of $p$ on $\Cal Q$ and let $B(\hat q, \epsilon)$ be a
closed $\epsilon$-ball defined by the total variation metric,
centered at $\hat q$. Then, for $\epsilon
> 0$,
 $$
\lim_{n \rightarrow \infty} \pi(q^n \in {B}(\hat q, \epsilon) \, |
\, (q^n \in \Cal Q) \wedge \nu^n) = 1.
 $$
\endproclaim

Proofs can be constructed along the lines for the static case.

 \subhead 4.2 General  prior
 \endsubhead

Let $\pi(q)$ be a prior pmf on $\Cal R$. 
From this pmf, a prior distribution $\pi^\Cal{A}(q^n)$ on
$\Cal{R}_n$ is constructed by a quantization $\Cal A \triangleq
\{A_1, A_2, \dots, A_N\}$ of $\Cal R$ into disjoint sets, such that
each $A \in \Cal A$ contains just one $q_n$ from $\Cal{R}_n$. Then
$\pi^{\Cal A}(q^n) \triangleq \pi(\{A_j: q^n \in A_j, j = 1, 2,
\dots, N\})$.

Let $\Cal S \triangleq S(\pi(\cdot))$. Let $\Cal{Q}^\pi \triangleq
\Cal Q \cap \Cal S$, $\Cal{P}^\pi \triangleq \Cal P \cap \Cal S$.

As the static case is subsumed under the dynamic one, only the
latter  limit theorems will be presented.

\proclaim {General prior $L$ST}
 Let $\nu^n \rightarrow p$. Let $\Cal Q$ be an open
set of sources. Let $\Cal{Q}^\pi \neq \emptyset$. Then,
 $$
\lim_{n \rightarrow \infty} \frac{1}{n} \log \pi^\Cal{A}(q^n \in Q|
\nu^n) = L(\Cal{Q}^\pi||p) - L(\Cal{P}^\pi||p).
 $$
\endproclaim

\demo {Proof} For a zero-prior-probability $n$-source, the posterior
probability is zero as well; so such sources can be excluded from
considerations.  Let $\Cal{S}_n \triangleq S(\pi^\Cal{A}(q_n))$,
$\Cal{Q}^\pi_n \triangleq \Cal{Q} \cap \Cal{S}_n$, $\Cal{P}^\pi_n
\triangleq \Cal{P} \cap \Cal{S}_n$.
$$
\log\pi^\Cal{A}(q^n \in \Cal Q| \nu^n) = 
\log{\sum_{q^n \in \Cal{Q}^\pi_n} \pi^\Cal{A}(q^n)\prod_\Cal{X}
(q^n)^{n \nu^n}} - \log{\sum_{q^n \in \Cal{P}^\pi_n}
\pi^\Cal{A}(q^n)\prod_\Cal{X} (q^n)^{n\nu^n}}.
$$

Denote by $\lambda(\Cal{Q}^\pi_n||\nu^n) \triangleq \sup_{q^n \in
\Cal{Q}^\pi_n} \lambda(q^n||\nu^n)$, where
$\lambda(q^n||\nu^n) 
\triangleq L(q^n||\nu^n) + \frac{1}{n}\log \pi^\Cal{A}(q^n)$. Using
this notation and invoking the same argument as in the proof of
$L$ST for uniform prior, $\frac{1}{n} \log \pi^\Cal{A}(q^{n} \in
\Cal{Q}| \nu^{n})$ can be bounded from above and below as:
$$\multline
\lambda(\Cal{Q}^\pi_n||\nu^{n}) - \lambda(\Cal{P}^\pi_n||\nu^n) -
\frac{m}{n}
\log(n + 1) \le \frac{1}{n} \log \pi^\Cal{A}(q^{n} \in \Cal Q| \nu^n) \le \\
\le \lambda(\Cal{Q}^\pi_n||\nu^{n}) -
\lambda(\Cal{P}^\pi_n||\nu^{n}) + \frac{m}{n} \log(n + 1).
\endmultline
$$
Since for $n\rightarrow\infty$, $\Cal{S}_n = \Cal S$, and $\nu^n
\rightarrow p$, and $\Cal Q$ is open, it taken together, implies
that $\lambda(\Cal{Q}^\pi_n||\nu^n)$ converges to
$L(\Cal{Q}^\pi||p)$. Similarly, $\lambda(\Cal{P}^\pi_n||\nu^n)$
converges to $L(\Cal{P}^\pi||p)$. $\square$
\enddemo

From the General prior $L$ST, follows
 \proclaim {General prior $L$CoLT}
Let $\nu^n \rightarrow p$. Let $\Cal Q$ be a convex, closed rare
(i.e., $p \notin \Cal Q$) set of sources. Let $\Cal{Q}^\pi \neq
\emptyset$ and let $\hat{q}^\pi$ be the $L$-projection of $p$ on
$\Cal{Q}^\pi$. Let $B(\hat{q}^\pi, \epsilon)$ be a closed
$\epsilon$-ball defined by the total variation metric, centered at
$\hat{q}^\pi$. Then, for $\epsilon
> 0$,
 $$
\lim_{n \rightarrow \infty} \pi^\Cal{A}(q^n \in {B}(\hat{q}^\pi,
\epsilon) \, | \, (q^n \in \Cal Q) \wedge \nu^n) = 1.
 $$
 \endproclaim

\head 4.3 Conditioning by rare sources: continuous alphabet
\endhead

Sanov's Theorem for continuous alphabet can be established either
via 'the method of types $+$ discrete approximation' approach (cf.
\cite{4}) or by means of the large deviations theory (cf. \cite{7}).
The former approach will be used here to formulate continuous
alphabet version of $L$ST.

Let $(\Cal Y, \Cal F)$ be a measurable space. Let $\Cal{T}^m$ be a
partition of the alphabet $\Cal Y$ into finite number $m$ of sets
$\Cal{T}^m \triangleq (T_1, T_2, \dots, T_m)$; $\Cal{T}_i \in \Cal
F$. The $\Cal{T}^m$-quantized $P$, denoted by $P^\Cal{T}$, is
defined as the distribution $P(T_1), P(T_2), \dots,$ $P(T_m)$ on the
finite set $\Cal{X} \triangleq \{1, 2, \dots, m\}$.

Let $\Cal P(\Cal Y)$ be the set of all probability measures on
$(\Cal Y, \Cal F)$. Let $\Cal Q \subseteq \Cal P$.  For probability
measures (pm's) $P, Q \in \Cal P(\Cal Y)$, the $L^m$-divergence
$L^m(Q||P)$ of $Q$ with respect to  $P$ is defined as
$$
 L^m(Q||P) \triangleq \sup_{\Cal{T}^m} L(Q^\Cal{T}||P^\Cal{T}),
$$ where the supremum is taken over all $m$-element partitions.
$L^m(\Cal{Q}||P)$ denotes  $\sup_{Q \in \Cal Q} L^m(Q||P)$. Let
$\Cal{Q}^\Cal{T} \triangleq \{Q: Q^\Cal{T} \in \Cal Q\}$,
$L^m(\Cal{Q}^\Cal{T}||P^\Cal{T}) \triangleq \sup_{
\Cal{Q}} L(Q^\Cal{T}||P^\Cal{T})$.

The empirical distribution $\nu^{n,m}$ of an $n$-sequence of $\Cal
Y$-valued random variables $Y$  with respect to a partition
$\Cal{T}^m$ is defined as  
$$
\nu^{n,m}_j = \frac{1}{n} \operatorname{Card}\{Y_i: Y_i \in T_j; \,
1\le i\le n\}, \qquad 1\le j \le m.
$$

The $\tau^m$-topology of pm's on $(\Cal Y, \Cal F)$ is the
topology in which a pm belongs to the interior of a set $\Cal Q$
of pm's iff for some partition $\Cal{T}^m$ and $\epsilon > 0$
$$
\{Q': |Q'(T_j) - Q(T_j)| < \epsilon, j = 1, 2, \dots, m\} \subset
\Cal Q.
$$

Thus, an $n$-source $q^n \in \Cal{R}_n(\Cal X)$
 belongs to the interior   of
$\Cal Q$ if there exists $\Cal{T}^m$ of $\Cal Y$ and $\epsilon
> 0$ such that the set $\{Q': |Q'({T}_j) - q^n_j| < \epsilon, j = 1,
2, \dots, m\}$ is a subset of $\Cal Q$.

Under the assumption of uniform prior distribution of $n$-sources, a
continuous analogue to the Dynamic $L$ST is:

\proclaim {Continuous $L$ST} Let, as $n \rightarrow\infty$,
$\nu^{n,m} \rightarrow R$, $R \in \Cal{R}(\Cal X)$. Let $\Cal Q$
be a rare (i.e., $R \notin \Cal Q$) open subset of $\Cal P(\Cal
Y)$. Then
$$
\lim_{n\rightarrow\infty} \frac{1}{n} \log \pi(q^n \in \Cal
Q|\nu^{n,m}) = L^m(\Cal{Q}||R) - L^m(\Cal{P}||R).
$$
\endproclaim

\demo {Proof}
First, an asymptotic lower bound to $\frac{1}{n}\log\pi(q^n \in \Cal
Q|\nu^n)$ will be established. Pick up a $Q$ such that for a
$\Cal{T}^m$, and an $\epsilon > 0$, $Q \in \Cal Q$. Let
$\Cal{M}^\Cal{T}(Q) \triangleq \{q^n: |q^n_j - Q(T_j)| < \epsilon, j
= 1, 2, \dots, m\}$. By the Dynamic $L$ST for uniform prior,
$\lim_{n\rightarrow\infty} \frac{1}{n}\log\pi(q^n \in
\Cal{M}^\Cal{T}(Q)|\nu^n) = L(\Cal{M}^\Cal{T}(Q)|R^\Cal{T}) -
L(R^\Cal{T}|R^\Cal{T})$ which is greater or equal to
$L(Q^\Cal{T}|R^\Cal{T}) - L(R^\Cal{T}|R^\Cal{T})$, since $Q^\Cal{T}
\in \Cal{M}^\Cal{T}(Q)$. Let $\Cal{M}(Q) \triangleq \cup_{\Cal{T}^m}
\Cal{M}^\Cal{T}(Q)$. Then, for $n \rightarrow \infty$,
$\frac{1}{n}\log\pi(q^n \in \Cal{M}(Q)|\nu^n) \ge \sup_{\Cal{T}^m}
L(Q^\Cal{T}|R^\Cal{T}) - L(R^\Cal{T}|R^\Cal{T}) \equiv L^m(Q||R) -
L^m(R||R)$. Since $\pi(q^n \in \Cal Q|\nu^n) \ge \sup_{Q \in
\Cal{Q}} \pi(q^n \in \Cal{M}(Q)|\nu^n)$,
$$
\lim_{n\rightarrow\infty}\frac{1}{n}\log\pi(q^n \in \Cal Q|\nu^n)
\ge \sup_{Q \in \Cal{Q}} L^m(Q||R) - L^m(R||R) \equiv
L^m(\Cal{Q}||R) - L^m(\Cal{P}||R).
$$

Asymptotic upper bound: for  $\Cal{T}^m$ as above, by the Dynamic
$L$ST with a uniform prior,
$$
\multline
 \lim_{n\rightarrow\infty}\frac{1}{n}\log\pi(q^n \in
\Cal{Q}^\Cal{T}|\nu^n) = L^m(\Cal{Q}^\Cal{T}||R^\Cal{T}) -
L^m(\Cal{P}^\Cal{T}||R^\Cal{T}) \\ \equiv \sup_{\Cal{Q}}
L(Q^\Cal{T}||\Cal{P}^\Cal{T}) - L(R^\Cal{T}||R^\Cal{T}).
 \endmultline
$$
Since $\pi(q^n \in \Cal{Q}|\nu^n) \le \sup_{\Cal{T}^m} \pi(q^n \in
\Cal{Q}^\Cal{T}|\nu^n)$,
$$
\lim_{n\rightarrow\infty}\frac{1}{n}\log\pi(q^n \in \Cal{Q}|\nu^n)
\le L^m(\Cal{Q}||R) - L^m(\Cal{P}||R).
$$
As the asymptotic lower and upper bounds coincide, the claim
follows. \qed
\enddemo

\head 5. Application to  Criterion Choice Problem
\endhead

1. Let there be an alphabet $\Cal X$ (finite, for simplicity) and
prior distribution $\pi(q^n)$ of $n$-rational sources. From the
prior $\pi(q^n$) an $n$-source is drawn, and the source then
generates an $n$-type $\nu^n$. We are not given the actual
$n$-source, but rather a set $\Cal Q$ to which the $n$-source
belongs. Given the alphabet $\Cal X$, the $n$-type $\nu^n$, the
prior distribution of sources $\pi(\cdot)$ and the set $\Cal Q
\subseteq \Cal P(\Cal X)$ the objective is to select an $n$-source
$q^n \in \Cal Q$. This constitutes the $\beta$-problem. Since
$\Cal Q$ in general contains more than one $n$-source the problem
is under-determined and in this sense ill-posed.

If $\Cal Q \equiv \Cal P(\Cal X)$, then under the assumption of
uniform prior distribution of $n$-sources, Static $L$LLN shows that
asymptotically (along the types $k$-equivalent with $\nu^n$) it is
just $\hat q \equiv \nu^n$ which is the 'only-possible' source of
$\nu^n$ (i.e., of itself)\footnote{Note that in the case of
unrestricted $\Cal Q$, $\nu^n$ is known to be Non-parametric Maximum
Likelihood Estimator of the source. Here, $\nu^n$ is the Maximum
A-posteriori Probability source.}. Dynamic $L$LLN, assuming that
$\nu^n \rightarrow r$, implies that the $n$-sources concentrate on
the true source. However, they do not, if a general prior is
assumed, such that it puts zero probability on the true source. In
the dynamic case ($\nu^n \rightarrow r$) with general prior,
$n$-sources concentrate on the $L$-projection of $r$ on
$\Cal{P}^\pi$.

What if $\Cal Q$ does not contain $\nu^n$? How should  an $n$-source
be selected in this case?
One possibility is to select $q^n$ from $\Cal Q$ by 
minimization of a distance or a convex statistical distance
measure \cite{16} between $\nu^n$ and $\Cal{Q}_n$. In this way,
the original $\beta$-problem of selecting $q^n \in \Cal Q$ is
transformed into an associated Criterion Choice Problem (CCP).

If the rare $\Cal Q$ is convex and closed, Static $L$CoLT shows
that - at least for $n$ sufficiently large - the CCP associated
with this instance  of the $\beta$-problem should be solved by
maximization of the $L$-divergence over $\Cal Q$. A major
qualifier has to be added to this statement: it holds provided
that uniform prior distribution of
$n$-sources is assumed. 
If a general prior, strictly positive on the entire set of rational
sources is assumed, then the statement still holds. Prior matters
only if it is not strictly positive on the entire $\Cal R$. Then, it
is the $L$-projection of $\nu^n$ on $\Cal{Q}^\pi$ that should be
selected (recall the General prior $L$CoLT).

2. Confront the $\beta$-problem with the following $\alpha$-problem
(also known as Boltzmann Jaynes Inverse Problem): let there be a
source $q$ that emits letters from an alphabet $\Cal X$. From the
source $q$ an $n$-type was drawn. We are not given the actual
$n$-type, but rather a set $\Pi$ to which the $n$-type belongs.
Given the alphabet $\Cal X$, the source $q$ and the set $\Cal{\Pi}$
the objective is to select an $n$-type $\nu^n \in \Cal\Pi$.

The CCP associated with the $\alpha$-problem is solved by CoLT and
GCP provided that $\Pi$ is a convex, closed rare set. The Theorems
imply that at least for sufficiently large $n$, the $I$-projection
of $q$ on $\Pi$ should be selected.

\head 6. EndNotes
\endhead

1) The terminology and notation of this paper follow more or less
closely \cite{2}, \cite{4}, \cite{6}, \cite{7}. The brief survey of
Large Deviations Theorems for Empirical Measures (Sect. 3) draws
from the same sources. For evolution of the results see among others
\cite{1}, \cite{3}, \cite{7}, \cite{8}, \cite{11}, \cite{12},
\cite{15}, \cite{17}, \cite{18}, \cite{21}, \cite{22}, \cite{23},
\cite{24}, \cite{25}. The inequalities used in Sect. 4 belong to
standard tool kit of the Method of Types, cf. \cite{4}. In relation
to the Proposition of Sect. 4.1 see also \cite{9}. The continuous
case of conditioning by rare sources (Sect. 5) is built parallel
with \cite{11} and \cite{4}.

2) This work is motivated by \cite{10}, where a problem of selecting
between Empirical Likelihood and Maximum Entropy Empirical
Likelihood (cf. \cite{19}, \cite{20}) has been addressed on
probabilistic, rather than statistical, grounds. Further discussion,
relevant also to the CCP associated with the $\alpha$ and
$\beta$-problems, can be found there.

3) Any of the results presented here may be stated in terms of
reverse $I$-projections \cite{5}. For instance the right-hand side
of the General prior $L$ST could be equivalently expressed as
$-(I(p||\Cal{Q}^\pi) - I(p||\Cal{P}^\pi))$, where $I(p||C)
\triangleq \inf_{q \in C} I(p||q)$ is the value of the
$I$-divergence at a reverse $I$-projection of $p$ on $C$. The above
mentioned statistical considerations (and 4) below) served as a
motivation 
for stating the results in terms of the newly introduced
$L$-divergence, though the $L$-projection is formally identical with
the reverse $I$-projection, which is already in use in a parametric
context, cf. \cite{5}. The present work leaves open the issue
whether it is more advantageous to state the Theorems of
conditioning by rare sources in terms of the $L$-projection or in
terms of the reverse $I$-projection.

4) If $p$ is an $n$-type then the $L$-divergence is known  as
Kerridge's inaccuracy; cf.~\cite{13},~\cite{14}.

5) For any prior $\pi(\cdot)$, the $L$-projection $\hat{q}^\pi$ of
$p$ on $\Cal{Q}^\pi$ is the same as the source which has
asymptotically supremal over $\Cal{Q}^\pi$ value of the posterior
probability $\pi(q^n|\nu^n)$. In the case of uniform prior the
correspondence holds for any $n$.

\enddocument